\documentclass[12pt]{amsart}

\usepackage{amsmath,amssymb,amscd,amsfonts,verbatim}
\usepackage{paralist}

\usepackage[mathscr]{eucal}

\newtheorem{thm}{Theorem}[section]
\newtheorem{lem}[thm]{Lemma}

\newtheorem{rems}[thm]{Remarks}

\newtheorem{prop}[thm]{Proposition}
\newtheorem{cor}[thm]{Corollary}
\newtheorem{defn}[thm]{Definition}

\newtheorem{assu-nota}[thm]{Assumption--Notation}
\theoremstyle{remark}

\newcommand{\C}{\mathbb C}

\newcommand{\OO}{\mathcal{O}}

\numberwithin{equation}{section}

\title[Base components]{ The base components of the dualizing sheaf of a curve on a surface}

\author[K.~Konno]{Kazuhiro Konno}
\address{Department of Mathematics, Graduate School of Science, 
Osaka University, Toyonaka, Osaka 560-0043, JAPAN}
\email{konno@math.sci.osaka-u.ac.jp}

\author[M.~Mendes Lopes]{Margarida Mendes Lopes}
\thanks{The second named author is a member of the Center for Mathematical
Analysis, Geometry and Dynamical Systems (IST-UTL).
This work was partially supported by the Funda\c{c}\~ao para a Ci\^encia e a Tecnologia
through Program POCI 2010/FEDER and Grants-in-Aid for Scientific Research (B) 
(No.~16340008) by Japan Society for the Promotion of Science (JSPS)}
\address{Center for Mathematical Analysis, Geometry, and Dynamical Systems,\\
Departamento de  Matem\'atica\\
Instituto Superior T\'ecnico\\
Universidade T{\'e}cnica de Lisboa\\
Av.~Rovisco Pais\\
1049-001 Lisboa, PORTUGAL}
\email{mmlopes@math.ist.utl.pt}



\begin{document}
\begin{abstract}
This note studies the  structure of the divisorial fixed part  of $|\omega_D|$  for a 1-connected curve $D$ on a smooth surface $S$.  
It is shown that if the divisorial fixed part $F$  of $|\omega_D|$ is non empty 
then it has arithmetic genus $\leq 0$ and each component of $F$ is a smooth rational curve.  
The stucture of curves $D$, with non empty divisorial fixed part $F$  for $|\omega_D|$,  is also described. 

\end{abstract}
\maketitle
\baselineskip=15.5pt

\section*{Introduction}

In this note we study the  structure of the fixed part  of $|\omega_D|$  for a 1-connected curve $D$ on a smooth surface $S$.  It is well known  that if $|\omega_D|$ has base points  then $D$ is not 2-connected (cf.~\cite{cfm}), 
and in fact there has been work of several authors concerning the structure of $\omega_D$ 
(see, e.g., \cite{cf}, \cite{cfm}, \cite{k},  \cite{m})  but  as far as we know the present result is new.

We prove the following theorem:

\begin{thm}\label{main} Let $D$ be a 1-connected curve on a smooth algebraic projective surface $S$ over $\C$. 
If $0\prec Z\prec D$ is a curve contained in the fixed part of  $|\omega_D|$ then:
\begin{enumerate}[\rm i)]

\item  every irreducible component of $Z$ is a smooth rational curve;
\item for any $0\prec Z'\preceq  Z$,  $p_a(Z')\leq 0$ and $h^1(Z',\OO_{Z'})=0$;
\item  $p_a(Z)=0$ if and only if $Z$ is 1-connected;
\item $h^0(D-Z,\OO_{D-Z})=(D-Z)Z+p_a(Z)$.
\end{enumerate}

Furthermore  if $p_a(Z)=0$ and $Z(D-Z)=m$, $D-Z$ decomposes as $B_1+ ...+B_m$, such that 
$\OO_{B_i}(-(B_{i+1}+... +B_m))=\OO_{B_i} $, for every $i=1, \dots ,m-1$,   $B_iZ=1$ and $B_i$ is 1-connected, for every $i=1, \dots, m$. In addition either $B_i\preceq B_{i+1}+...+B_m$ or $B_i\cap B_{i+1}+ ...+B_m=\emptyset$.
\end{thm}

\bigskip

\section{Preliminaries}
\subsection  {Notation.} 

By a  {\it curve} we mean a non-zero effective divisor on a smooth algebraic projective surface $S$ over $\C$.  $K$ will denote the canonical bundle of $S$.

Given a curve $D$, $\omega_D$ denotes its dualizing sheaf and $p_a(D)$ its arithmetic genus.

 A curve  $D$ is $m$-connected if for every decomposition $D=A+B$, with $A,B$ effective non-zero divisors, $AB \geq m$.
 
 For any invertible sheaf  $\mathcal{L}$ on $D$ we denote by 
$h^i(D,\mathcal{L})$ the  dimension as a $\C$-vector space  of $H^i(D,\mathcal{L})$.

\subsection{Some properties} 
 Here we list some properties that will be used throughout without further reference.

\begin{itemize} 
\item Given a curve  $D$, $2p_a(D)-2= KD+D^2$ (adjunction formula).
\item $h^0(D,\OO_D)-h^1(D,\OO_D)=1-p_a(D)$. %where $2p_a(D)-2=KD+D^2$.
\item By duality one has $h^0(D,\OO_D)=h^1(D,\omega_D)$ and $h^1(D,\OO_D)=h^0(D,\omega_D)$. 
\item If a curve $D$ decomposes as the sum of two curves $D_1,D_2$ then
$2p_a(D)-2=KD+D^2=KD_1+KD_2+D_1^2+D_2^2+2D_1D_2=2p_a(D_1)-2+2p_a(D_2)-2+2D_1D_2$
and so $p_a(D)=p_a(D_1)+p_a(D_2)-1+D_1D_2$.
\item (\cite{cfm}) Let $D$ be an
$m$-connected curve on the surface $S$ and let $D=D_1+D_2$ with $D_1$, $D_2$
curves. Then:\begin{enumerate}[\rm i)]
\item if  $D_1 \cdot D_2=m$, then $D_1$ and $D_2$ are
$[(m+1)/2]$-connected;
\item if $D_1$ is chosen to be minimal subject to the
condition $D_1 \cdot (D-D_1) = m$, then $D_1$ is $[(m+3)/2]$-connected.
\end{enumerate} 
\end{itemize}

\subsection{Auxiliary results} 
\begin{defn} Let $D$ be a reducible
curve on a smooth surface $S$, ${\mathcal L}$ an invertible sheaf on $D$ and let
$s\in H^0(D,{\mathcal L})$ with $s\neq 0$ such that $s$ vanishes identically on some
component of $D$. Let  $Z_s\prec D$ be the
biggest curve such that $s_{|Z_s}\equiv 0$.\par We will say that $s$ is
{\it 0-maximal} if there is no global section $t$ of ${\mathcal L}$ with $Z_s\prec Z_t$.
\end{defn}

\begin{lem}\label{dec} Let $A$ be a curve  such that $h^0(A,\OO_A)\geq 2$. Then there is a decomposition $A=A_1+A_2$ where $A_1$, $A_2$ are curves such that
\begin{enumerate}[\rm i)]

\item  $h^0(A_1,\OO_{A_1}(-A_2))\neq 0$;

\item $A_1A_2\leq 0$;
\item for each component $\Gamma$ of $A_1$, $\Gamma A_2\leq 0$;
\item for each component $\Gamma$ of $A_1$, the restriction map $$H^0(A_1,\OO_{A_1}(-A_2))\to H^0(\Gamma,\OO_{\Gamma}(-A_2))$$ is injective;

\item  for each component $\Gamma$ of $A_1$, $$h^0(A,\OO_A)\leq h^0(A_2,\OO_{A_2})+h^0(\Gamma,\OO_{\Gamma}(-A_2)).$$

 \end{enumerate}
\smallskip 
Furthermore if $A_1A_2=0$ then $\OO_{A_1}(-A_2)=\OO_{A_1}$, and $h^0(A_1,\OO_{A_1})=1$.
\end{lem}

\begin{proof}

Since  $h^0(A,\OO_A) \geq 2 $,
there exists a section $s$ in $H^0(A, \OO_A)$ vanishing identically on some
component of $A$. Choose a 0-maximal such section $s$ and let $A_2:=Z_s$,
$A_1:=D-A_2$. From the exact sequence\par
 $$0 \to \OO_{A_1}(-A_2) \to \OO_A\to \OO_{A_2} \to 0$$ \noindent we get\par 

$$0 \to H^0(A_1,\OO_{A_1}(-A_2)
) \to H^0(A,\OO_A)\buildrel \ r \over \to
H^0(A_2,\OO_{A_2}). $$\par

By the hypothesis of 0-maximality of $s$,  every section of $ H^0(A_1,\OO_{A_1}(-A_2))$ does not vanish identically on any component $\Gamma$ of $A_1$ and so in particular  $\Gamma A_2\leq 0$ and $A_1A_2\leq 0$.  This proves assertions   ii) and iii). 

Furthermore for any $\Gamma$, again the hypothesis of 0-maximality implies that the kernel of the restriction map   
$$H^0(A_1, \OO_{A_1}(-A_2))\to  H^0(\Gamma, \OO_{\Gamma}(-A_2))$$ is 
0-dimensional and therefore we get  assertions  iv) and v), because 

$$h^0(A,\OO_A)\leq h^0(A_1,\OO_{A_1}(-A_2))+h^0(A_2,\OO_{A_2}).$$ 

The last assertion is clear, by the previous considerations.
\end{proof}

\medskip

\begin{rems} 

\par\noindent{\rm  a)  In the decomposition above, if there is a component  $\Gamma$ of $A_1$ such that $\Gamma A_2=0$, then $ h^0(A_1,\OO_{A_1}(-A_2))= 1$.  Otherwise again one  would get a contradiction to $0$-maximality of $s$.

\noindent  b)   Lemma \ref{dec} means that a curve $D$ such that  $h^0(D,\OO_D) \geq 2 $  is  necessarily not 1-connected. 
However  $h^0(D,\OO_D)=1$ does not imply 1-connectedness. 
For instance a multiple fibre  $F=mD$ of a fibration is not 1-connected  but  $h^0(F,\OO_F)=1$ (cf.~\cite[Chp. III]{bpv} )}
\end{rems}

\medskip

\begin{lem}\label{b} Let  $D$ be a 1-connected curve. If $A\prec D$ is such that $A(D-A)=b$, then $h^0(A,\OO_A)\leq b$.
\end{lem}
\begin{proof}
We do this proof by induction on $b$.
If $b=1$ then, because $D$ is 1-connected, $A$ is also and so $h^0(A,\OO_A)=1$.
We assume that we have proved the assertion for $m<b$ and we want to prove for $m=b$.

Suppose that $h^0(A,\OO_A)\geq 2$.
Take a decomposition of $A$ as in  Lemma \ref{dec} and let $A_1A_2=-\alpha$, with $\alpha\geq 0$. Note that, because for  every component $\Gamma$ of $A_1$, $\Gamma A_2\leq 0$, we have, for any component $\Gamma$ of $A_1$, $\Gamma A_2\geq -\alpha$, i.e. $\Gamma(-A_2)\leq \alpha$, and so $h^0(\Gamma, \OO_{\Gamma}(-A_2))\leq \alpha+1$. 

Now, by 1-connectedness of $D$, $A_1(A_2+(D-A))\geq 1$, and so $A_1(D-A)\geq 1+\alpha$. Since $(A_1+A_2)(D-A)=b$, $A_2(D-A)\leq b-1-\alpha$ and so $A_2(A_1+(D-A))=A_2(D-A_2)\leq b-1-2\alpha$. Hence by the induction hypothesis 
$h^0(A_2,\OO_{A_2})\leq b-1-2\alpha$ and so by  Lemma \ref{dec}
$h^0(A,\OO_A)\leq b-\alpha$.
\end{proof}

\smallskip

\begin{cor} \label{h0} Suppose that  the curve $D$ is 1-connected. 
If $A\prec D$ is such that $A(D-A)=b$ and $h^0(A,\OO_A)= b$, then $A$ decomposes as $B_1+....+B_b$, 
such that $\OO_{B_i}(-(B_{i+1}+....+B_b))=\OO_{B_i} $,  for every $i=1,\dots ,b-1$, $B_i(D-A)=1$ and $B_i$ is 1-connected, for every $i=1, \dots ,b$. 
Furthermore either $B_i\preceq B_{i+1}+....+B_b$ or $B_i\cap B_{i+1}+....+B_b=\emptyset$.
\end{cor}

\begin{proof}
Suppose that $h^0(A,\OO_A)= b$. 
Then in the decomposition as in the previous proof we must have $\alpha=0$ and $A_1(D-A)=1$, meaning that $A_2(D-A)=b-1$. 
So by Lemma \ref{dec} necessarily $h^0(A_1(-A_2)=1$ and $\OO_{A_1}(-A_2)=\OO_{A_1}$.

Note also that $A_1(D-A_1)=1$ means that $A_1$ is 1-connected. Assume that 
$A_1$ has common components with $A_2$. 
Then we can write $A_1=H+B, A_2=H+C$ where $B$ and $C$ have no common components and  $H\neq 0$. 
Suppose $B\neq 0$. 
Since $\OO_{A_1}(-A_2)=\OO_{A_1}$, $B(H+C)=0$ and so, because $BC\geq 0$, we conclude that $BH\leq 0$. 
But this contradicts the 1-connectedness of $A_1$ and so $B=0$.

We take $B_1:=A_1$. Now we consider $D-A_1$ which is still 1-connected. 
One has $h^0(A_2,\OO_{A_2})=b-1$ and $A_2(D-A_1-A_2)=b-1$. 
We can apply the same reasoning as before and an obvious induction gives us the statement.
\end{proof}

\smallskip

\begin{lem}\label{pa1}
Let \(\mathcal{L}\) be an invertible sheaf on a curve \(Z\) satisfying 
\(\deg \mathcal{L}|_{Z'}\geq 2p_a(Z')-2\) for any \(0\prec Z'\preceq Z\).
If \(h^1(Z,\mathcal{L})\neq 0\), then there exists a subcurve \(A \preceq Z\) 
with \(\mathcal{L}|_A\simeq 
 \omega_{A}\) and \(h^0(A,\mathcal{O}_{A})=1\).
\end{lem}

\begin{proof}
By duality, we have \(h^0(Z,\omega_{Z}\otimes \mathcal{L}^{-1})\neq 0\).
Take a \(0\)-maximal \(s\in H^0(Z,\omega_Z \otimes \mathcal{L}^{-1})\) and 
put \(A=Z-Z_s\) (possibly \(Z_s=0\)). 
Then \(\mathcal{O}_{A}(-Z_s)\otimes \omega_Z\otimes \mathcal{L}^{-1}\simeq 
\omega_A\otimes \mathcal{L}^{-1}\) and \(s\) induces a non-zero 
\(s'\in H^0(A,\omega_A\otimes \mathcal{L}^{-1})\) which does 
 not vanish identically on any component of \(A\).
In particular, \(\omega_A\otimes \mathcal{L}^{-1}\) is nef.
Since \(\deg \omega_A=2p_a(A)-2\leq \deg \mathcal{L}|_{A}\) by the assumption, 
\(\omega_A\otimes \mathcal{L}^{-1}\) is numerically trivial.
Furthermore, since \(s'\) is nowhere vanishing, 
we get \(\omega_A\otimes \mathcal{L}^{-1}\simeq \mathcal{O}_{A}\).
By the \(0\)-maximality of \(s\), the restriction map 
 \(H^0(A,\omega_A\otimes \mathcal{L}^{-1})=H^0(A,\mathcal{O}_{A})\to H^0(\Gamma,\mathcal{O}_{\Gamma})\) 
is injective for any irreducible component \(\Gamma \preceq A\).
Hence \(h^0(A,\mathcal{O}_{A})=1\).
\end{proof}

\smallskip

\begin{prop} \label{go} 
Let $Z$ be a curve such that for any $0\prec Z'\preceq  Z$, $p_a(Z')\leq 0$. 
Then:
\begin{enumerate}[\rm i)]

\item  every component $\Gamma$ of $Z$ is a smooth rational curve;

\item for any $0\prec Z'\preceq  Z$, $h^1(Z',\OO_{Z'})=0$.
\end{enumerate}
Furthermore, \(p_a(Z)=0\)  (and $h^0(Z,\OO_Z)=1$)  if and only if \(Z\) is 1-connected.
\end{prop}

\begin{proof}
By Lemma~\ref{pa1} applied to \(\mathcal{L}=\mathcal{O}_Z\), if 
\(h^1(Z,\mathcal{O}_Z)\neq 0\) there is  a subcurve \(A\) with 
\(\mathcal{O}_A\simeq \omega_A\) and \(h^0(A,\mathcal{O}_A)=1\).
This implies \(p_a(A)=1\) contradicting the hypothesis \(p_a(A)\leq 0\).
Therefore, \(h^1(Z,\mathcal{O}_Z)=0\).
Then we get i) and ii), since the natural map \(H^1(Z,\mathcal{O}_Z)\to 
H^1(Z',\mathcal{O}_{Z'})\) is surjective for any \(0\prec Z'\preceq Z\).

If \(Z\) is \(1\)-connected, then \(h^0(Z,\mathcal{O}_Z)=1\)  (cf. Lemma \ref{dec}) and we have
\(p_a(Z)=0\) by \(h^1(Z,\mathcal{O}_Z)=0\).
Conversely, assume that \(p_a(Z)=0\). 
Then \(0=p_a(Z)=p_a(Z_1)+p_a(Z_2)-1+Z_1Z_2\leq -1+Z_1Z_2\) 
for any decomposition \(Z=Z_1+Z_2\) with \(0\prec Z_1,Z_2\).
Hence \(Z\) is \(1\)-connected.
\end{proof}

\bigskip

\section{Fixed components of $|\omega_D|$}

The results in this section prove Theorem \ref{main}.  

\smallskip

\begin{prop}\label{pa} Let $D$ be a 1-connected curve  and $Z\prec D$ a  curve such that the restriction map $H^0(D, \omega_D)\to H^0(Z, \omega_D)$ is the zero map.
Then $p_a(Z)\leq 0$ and $h^0(D-Z,\OO_{D-Z})=(D-Z)Z+p_a(Z)$.
\end{prop}

\begin{proof}
Let $B:=D-Z$. 
As usual one has $p_a(D)= p_a(B)+p_a(Z)-1+BZ$, and so $p_a(B)-1=p_a(D)-BZ-p_a(Z)$.

Since the kernel of the restriction map $H^0(D, \omega_D)\to H^0(Z, \omega_D)$ 
is exactly $H^0(B, \omega_B)$, our hypothesis implies that $h^0(B, \omega_B)=h^0(D, \omega_D)$. 
The equality  $p_a(B)-1=h^0(B, \omega_B)-h^0(B,\OO_B)$ yields then 
$$p_a(D)-h^0(B,\OO_B)=p_a(D)-BZ-p_a(Z),$$\noindent i.e., 
$h^0(B,\OO_B)=BZ+p_a(Z)$.

By Lemma \ref{b}, $h^0(B,\OO_B)\leq BZ$ and so $p_a(Z)\leq 0$ .
\end{proof}

\smallskip

\begin{cor}\label{h1}  
Let $D$ be a 1-connected curve and let $Z$ be the fixed part of $|\omega_D|$. 
Then for any $Z'\preceq  Z$, $p_a(Z')\leq 0$ and $h^1(Z', \OO_{Z'})=0$.
\end{cor}
\begin{proof}
The statement is an immediate consequence of Propositions \ref{pa} and  \ref{go}.
\end{proof}

\smallskip

\begin{cor} 
Let $D$ be a 1-connected curve  and $Z\prec D$ a 1-connected  curve (or more generally such that $h^0(Z,\OO_Z)=1$) 
such that the restriction map $H^0(D, \omega_D)\to H^0(Z, \omega_D)$ is the zero map.
Then $h^1(Z,\OO_Z)= 0$ and $B:=D-Z$ decomposes as in Corollary \ref{h0}.
\end{cor}

\begin{proof}  By Corollary \ref{h1}, $p_a(Z)\leq 0$ and    $h^1(Z,\OO_Z)= 0$,  and so, because  $h^0(Z,\OO_Z)=1$, 
we have $p_a(Z)=0$. 
By Proposition \ref{pa},  $h^0(B,\OO_B)=BZ$ and therefore we can apply Corollary \ref{h0} obtaining  a decomposition as wished.
\end{proof}

%\noindent{\bf Acknowledgements.} %I want to give special thanks to Kazuhiro Konno who made me think again of these problems.
% This work was partially supported by the Funda\c{c}\~ao para a Ci\^encia e a Tecnologia
%through Program POCI 2010/FEDER.

\bigskip

%\begin{minipage}{13cm}
%\parbox[t]{6.5cm}{Margarida Mendes Lopes\\
%Center for Mathematical Analysis, Geometry, and Dynamical Systems,\\
%Departamento de  Matem\'atica\\
%Instituto Superior T\'ecnico\\
%Universidade T{\'e}cnica de Lisboa\\
%Av.~Rovisco Pais\\
%1049-001 Lisboa, PORTUGAL\\
%mmlopes@math.ist.utl.pt
%  } \hfill

%\end{minipage}


\bigskip

\begin{thebibliography}{ABCD}
\bibitem[BPV]{bpv}  W.~Barth, C.~Peters, A.~Van De Ven, {\em Compact complex surfaces},
Ergebnisse der Mathematik und ihrer Grenzgebiete, {\bf3.} Folge, Band {\bf 4},
Springer-Verlag, Berlin (1984).

\bibitem[CCFR]{cf} 
F.~Catanese, M.~Franciosi, K.~Hulek, M.~Reid, 
{\em Embeddings of curves and surfaces},  
Nagoya Math. J. {\bf 154}  (1999), 185--220. 

\bibitem[CFM]{cfm} 
C. Ciliberto, P. Francia, M. Mendes Lopes, 
{\em Remarks on the bicanonical map for surfaces of general type}, 
Math. Z. {\bf 224} (1997), 137--166.

\bibitem [K]{k} K.~Konno, 
{\em 1-2-3 theorem for curves on algebraic surfaces},  
J. Reine Angew. Math.  {\bf 533}  (2001), 171--205. 

\bibitem [M] {m} M.~Mendes Lopes, 
{\em Adjoint systems on surfaces}, 
Boll. Un. Mat. Ital. A (7)  {\bf 10}  (1996),  no. 1, 169--179.

\end{thebibliography}
\end{document}